\documentclass[a4paper,11pt]{scrartcl}
\usepackage{amsmath,amssymb}
\usepackage{amsthm}
\usepackage{graphicx}
\usepackage[square,sort&compress,numbers]{natbib}
\usepackage{url}
\usepackage{type1cm}
\usepackage{booktabs}

\usepackage{subcaption}
\usepackage{newtxtext,newtxmath}

\newcommand{\Authornote}{\renewcommand{\thefootnote}{\fnsymbol{footnote}}}
\newcommand{\authornote}{\Authornote\footnote}
\theoremstyle{plain}
\newtheorem{theorem}{Theorem}[section]

\theoremstyle{definition}

\theoremstyle{remark}
\newtheorem{remark}[theorem]{Remark}

\newcommand{\reffig}[1]{Figure~\ref{#1}}
\newcommand{\reftab}[1]{Table~\ref{#1}}

\newcommand{\finbox}{\nolinebreak\hfill{\small $\blacksquare$}}

\newcommand{\MIN}{\mathop{\mathrm{Minimize}}}

\newcommand{\ST}{\mathop{\mathrm{subject~to}}}

\renewcommand{\Re}{\ensuremath{\mathbb{R}}}

\newcommand{\bi}[1]{\ensuremath{\boldsymbol{#1}}}

\newcommand{\bs}[1]{\ensuremath{\boldsymbol{\mathsf{#1}}}}

\begin{document}

\begin{center}
  {\Large\bfseries\sffamily%
  Mixed-Integer Programming Formulation }\\
  \medskip
  {\Large\bfseries\sffamily%
  of a Data-Driven Solver in Computational Elasticity 
  }%
  \par%
  \bigskip%
  {
  Yoshihiro Kanno~\authornote[2]{%
    Mathematics and Informatics Center, 
    The University of Tokyo, 
    Hongo 7-3-1, Tokyo 113-8656, Japan.
    E-mail: \texttt{kanno@mist.i.u-tokyo.ac.jp}. 
    }
  }
\end{center}

\begin{abstract}
  This paper presents a mixed-integer quadratic programming formulation 
  of an existing data-driven approach to computational elasticity. 
  This formulation is suitable for application of a standard 
  mixed-integer programming solver, which finds a global optimal 
  solution. 
  Therefore, the results obtained by the presented method can be used as 
  benchmark instances for any other algorithm. 
  Preliminary numerical experiments are performed to compare quality of 
  solutions obtained by the proposed method and a heuristic used in the 
  data-driven computational mechanics. 
\end{abstract}

\begin{quote}
  \textbf{Keywords}
  \par
  Global optimization; 
  mixed-integer programming; 
  exact solution method; 
  data-driven computing; 
  model-free computational mechanics. 
\end{quote}

\section{Introduction}

Data-driven computation in elasticity was initiated by \citet{KO16}. 
This methodology directly uses the material data set obtained from 
physical experiments, rather resorting to conventional empirical 
modeling of a material constitutive law. 
It attempts to minimize the distance between the data set and 
the strains and the stresses satisfying the compatibility relation and 
the force-balance equation. 
This method has recently been extended to 
geometrically nonlinear problems \cite{NK18} 
and dynamic problems \cite{KO18}. 
Also, \citet{KO17} introduced the information entropy to reduce the 
variance of the original method in \cite{KO16}. 

The method in \cite{KO16} is regarded as a lazy learning method, where 
no model is learned from the given data set. 
In contrast, another data-driven approach proposed by 
\citet{IAcAGCC18,IBAAcCLC17} is based on the manifold learning, which is 
one of eager learning methods. 
Also, to reduce the influence of outliers in a material data set, a method 
using the local robust regression has been proposed for static analysis 
of trusses \cite{Kan18}. 

Attention of this note is focused on a numerical solution of the 
data-driven approach in \cite{KO16}. 
The problem dealt with in \cite{KO16} is essentially an optimization 
problem. 
Algorithm~1 in \cite{KO16} serves as a heuristic for this optimization problem. 
Subsequently, pretty much the same heuristics have been used for 
data-driven computational mechanics \cite{NK18,LCRSV18}. 
In this note, we show that this optimization problem can be formulated 
as a {\em mixed-integer quadratic programming\/} (MIQP) problem in a 
natural manner. 
This problem can be solved {\em globally\/} with, 
e.g., a branch-and-bound method, 
because its continuous relaxation is a (convex) quadratic programming problem. 
Several sophisticated software packages are available for solving MIQP problems. 
Although the modeling presented in this note is fairly standard in 
integer optimization, it cannot be found in literature to the best of 
the author's knowledge.

\section{Problem setting}
\label{sec:problem}
In this section, we overview the methodology of data-driven computing in \cite{KO16}. 
Although the formulation presented in section~\ref{sec:mixed} can 
readily be adapted to more general structures, e.g., three-dimensional 
continua, we restrict ourselves to trusses for simple presentation. 
Throughout the paper, we assume elasticity and small deformation. 

Suppose that we are given experimental material data, consisting of 
pairs of uniaxial strain and uniaxial stress values. 
We use 
$D = \{ (\check{\varepsilon}_{1},\check{\sigma}_{1}),
  \dots, (\check{\varepsilon}_{d},\check{\sigma}_{d}) \}$ 
to denote the data set, 
where $\check{\varepsilon}_{j}$ and $\check{\sigma}_{j}$ are the 
observed strain and stress, respectively, 
and $d$ is the number of observations. 

Consider a truss consisting of this material. 
Let $m$ and $n$ denote the number of members and the number of degrees 
of freedom of the nodal displacements, respectively. 
We use $\varepsilon_{i} \in \Re$ and $\bi{u} \in \Re^{n}$ to denote 
the axial strain of member $i$ and the nodal displacement vector, respectively. 
The compatibility relations can be described in the form 
\begin{align}
  \varepsilon_{i} = \bi{b}_{i}^{\top} \bi{u} , 
  \quad  i=1,\dots,m , 
  \label{eq:constitutive}
\end{align}
where $\bi{b}_{i} \in \Re^{n}$ is a constant vector. 

Let $\sigma_{i} \in \Re$ and $\bi{p} \in \Re^{n}$ denote the axial 
stress of member $i$ and the external load vector, respectively. 
The force-balance equations are written as 
\begin{align}
  \sum_{i=1}^{m} v_{i} \sigma_{i} \bi{b}_{i} = \bi{p} , 
  \label{eq:force-balance}
\end{align}
where $v_{i}$ is the volume of member $i$. 

To state the methodology of data-driven computing in \cite{KO16} 
succinctly, we attempt to find the points 
$(\varepsilon_{1},\sigma_{1}),\dots,(\varepsilon_{m},\sigma_{m})$ that 
are ``closest'' to $D$, when \eqref{eq:constitutive} and 
\eqref{eq:force-balance} are satisfied. 
Define the distance from point $(\varepsilon_{i},\sigma_{i})$ to 
data set $D$ by 
\begin{align}
  f(\varepsilon_{i},\sigma_{i}) 
  = \min 
  \left\{\left.
  \sqrt{ \frac{v_{i}}{2} }
  \begin{Vmatrix}
    \begin{bmatrix}
      \sqrt{c} (\varepsilon_{i} - \check{\varepsilon}) \\
      (\sigma_{i} - \check{\sigma}) / \sqrt{c} \\
    \end{bmatrix}
  \end{Vmatrix}
  _{2}
  \ \right|
  (\check{\varepsilon}, \check{\sigma}) \in D 
  \right\} , 
  \label{eq:distance}
\end{align}
where $c > 0$ is a constant, and 
$\| \,\cdot\, \|_{2}$ denotes the Euclidean norm of a vector. 
Then we we minimize the sum of squared distances, i.e., 
$\sum_{i=1}^{m} f(\varepsilon_{i},\sigma_{i})^{2}$.

\section{Mixed-integer quadratic programming formulation}
\label{sec:mixed}

In this section, we recast the problem described in 
section~\ref{sec:problem} as an MIQP problem. 

Observe that, in the minimization problem in \eqref{eq:distance}, 
we select one data point to evaluate the distance from 
$(\varepsilon_{i},\sigma_{i})$ to $D$. 
We use 0-1 variables, $t_{i1},\dots,t_{id}$, to represent this selection 
such that $t_{ij}=1$ if 
$(\check{\varepsilon}_{j},\check{\sigma}_{j})$ is selected, and 
otherwise $t_{ij}=0$. 
We see that $(e_{i},s_{i}) \in D$ if and only if there exist 
$t_{i1},\dots,t_{id}$ satisfying 
\begin{align}
  \begin{bmatrix}
    e_{i} \\ s_{i} \\
  \end{bmatrix}
  &= \sum_{j=1}^{d} 
  \begin{bmatrix}
    \check{\varepsilon}_{j} \\
    \check{\sigma}_{j} \\
  \end{bmatrix}
  t_{ij} , 
  \label{eq:SOS.1} \\
  \sum_{j=1}^{d} t_{ij} &= 1 , 
  \label{eq:SOS.2} \\
  t_{ij} & \in \{ 0,1 \} , 
  \quad  j=1,\dots,d . 
  \label{eq:SOS.3}
\end{align}
From \eqref{eq:distance}, we obtain 
\begin{align}
  f(\varepsilon_{i},\sigma_{i})^{2} 
  = \min_{e_{i},s_{i},t_{i1},\dots,t_{id}} \left\{\left.
  \frac{1}{2} v_{i} c (\varepsilon_{i} - e_{i})^{2} 
  + \frac{1}{2} \frac{v_{i}}{c} (\sigma_{i} - s_{i})^{2} 
  \ \right|
  \eqref{eq:SOS.1}, \eqref{eq:SOS.2} , \eqref{eq:SOS.3}
  \right\} . 
  \label{eq:square.distance}
\end{align}
It is worth noting that, at the optimal solution of 
\eqref{eq:square.distance}, we have 
$t_{ij}=1$ if $(\check{\varepsilon}_{j},\check{\sigma}_{j})$ is 
the closest to $(\varepsilon_{i},\sigma_{i})$, 
and otherwise $t_{ij}=0$. 

In the data-driven solver \cite{KO16}, we minimize the sum of 
$f(\varepsilon_{i},\sigma_{i})^{2}$ in \eqref{eq:square.distance} under 
the constraints in \eqref{eq:constitutive} and \eqref{eq:force-balance}. 
This problem can be recast as follows: 
\begin{subequations}\label{P.original.1}%
  \begin{alignat}{3}
    & \MIN  &{\quad}& 
    \sum_{i=1}^{m} \frac{1}{2} v_{i} c (\varepsilon_{i} - e_{i})^{2} 
    + \sum_{i=1}^{m} \frac{1}{2} \frac{v_{i}}{c} (\sigma_{i} - s_{i})^{2} 
    \label{P.original.1.1} \\
    & \ST && 
    \begin{bmatrix}
      e_{i} \\ s_{i} \\
    \end{bmatrix}
    = \sum_{j=1}^{d} 
    \begin{bmatrix}
      \check{\varepsilon}_{j} \\
      \check{\sigma}_{j} \\
    \end{bmatrix}
    t_{ij} , 
    {\quad} i=1,\dots,m, 
    \label{P.original.1.2} \\
    & &&
    \sum_{j=1}^{d} t_{ij} = 1 , 
    {\quad} i=1,\dots,m, \\
    & &&
    t_{ij} \in \{ 0,1 \} , 
    {\quad} i=1,\dots,m; \ j=1,\dots,d , \\
    & &&
    \varepsilon_{i} = \bi{b}_{i}^{\top} \bi{u}, 
    {\quad} i=1,\dots,m, 
    \label{P.original.1.5} \\
    & &&
    \sum_{i=1}^{m} v_{i} \sigma_{i} \bi{b}_{i}= \bi{p} . 
    \label{P.original.1.6}
  \end{alignat}
\end{subequations}
Here, variables to be optimized are 
$\bi{u} \in \Re^{n}$, $\varepsilon_{i}$, $\sigma_{i}$, $e_{i}$, $s_{i}$, 
and $t_{ij}$ $(i=1,\dots,m;\, j=1,\dots,d)$. 

Problem \eqref{P.original.1} is minimization of a convex quadratic 
function under some linear equality constraints and some binary constraints. 
Hence, it is an MIQP problem. 
This problem can be solved globally with a standard MIP 
(mixed-integer programming) solver. 

\begin{remark}
  It is fairly straightforward to extend problem \eqref{P.original.1} to 
  a continuum discretized by the conventional finite element method. 
  The material data set is now 
  $D = \{ (\check{\bi{\varepsilon}}_{1},\check{\bi{\sigma}}_{1}),
  \dots, (\check{\bi{\varepsilon}}_{d},\check{\bi{\sigma}}_{d}) \}$, 
  where $\check{\bi{\varepsilon}}_{j}$ and $\check{\bi{\sigma}}_{j}$ are 
  second-order symmetric tensors with the dimension three. 
  Instead of the member strain and stress, we attempt to compute the 
  strain and stress tensors at each evaluation point of numerical 
  integration, as well as the nodal displacement vector. 
  The objective function in \eqref{P.original.1.1} is then replaced by 
  \begin{align*}
    \sum_{i=1}^{m} \frac{1}{2} \rho_{i} \bs{C} 
    (\bi{\varepsilon}_{i} - \bi{e}_{i}) : (\bi{\varepsilon}_{i} - \bi{e}_{i}) 
    + \sum_{i=1}^{m} \frac{1}{2} \rho_{i} \bs{C}^{-1} 
    (\bi{\sigma}_{i} - \bi{s}_{i}) : (\bi{\sigma}_{i} - \bi{s}_{i}) , 
  \end{align*}
  where $\bs{C}$ is a constant positive definite fourth-order tensor, 
  and $\rho_{i}>0$ is the weight for the numerical integration. 
  Thus, the objective function is still a convex quadratic function. 
  Also, the constitutive relations in \eqref{P.original.1.5} and the 
  force-balance equations in \eqref{P.original.1.6} remain to be linear 
  equality constraints, because we assume small deformation. 
  Obvious change in \eqref{P.original.1.2} keeps its linearity. 
  Thus, the problem for continua is also an MIQP problem. 
  \finbox
\end{remark}

\section{Numerical experiments}
\label{sec:ex}
\begin{figure}[tp]
  \centering
  \includegraphics[scale=0.90]{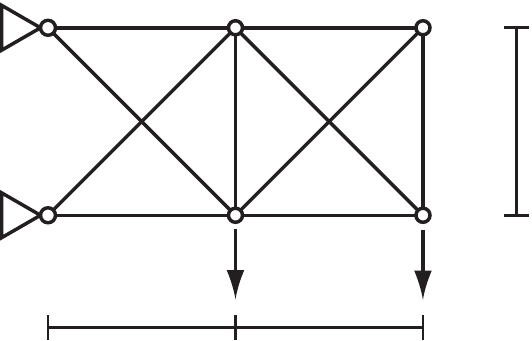}
  \begin{picture}(0,0)
    \put(-150,-50){
    \put(36,45){{\footnotesize $3.6\,\mathrm{m}$ }}
    \put(84,45){{\footnotesize $3.6\,\mathrm{m}$ }}
    \put(145,104){{\footnotesize $3.6\,\mathrm{m}$ }}
    }
  \end{picture}
  \medskip
  \caption{A 10-bar truss. }
  \label{fig:10bar}
\end{figure}
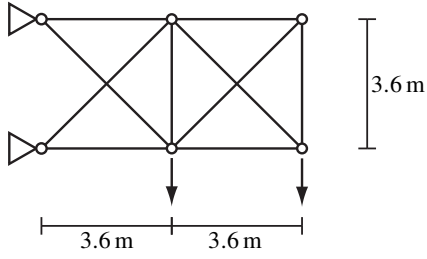

\begin{figure}[tp]
  \centering
  \includegraphics[scale=0.45]{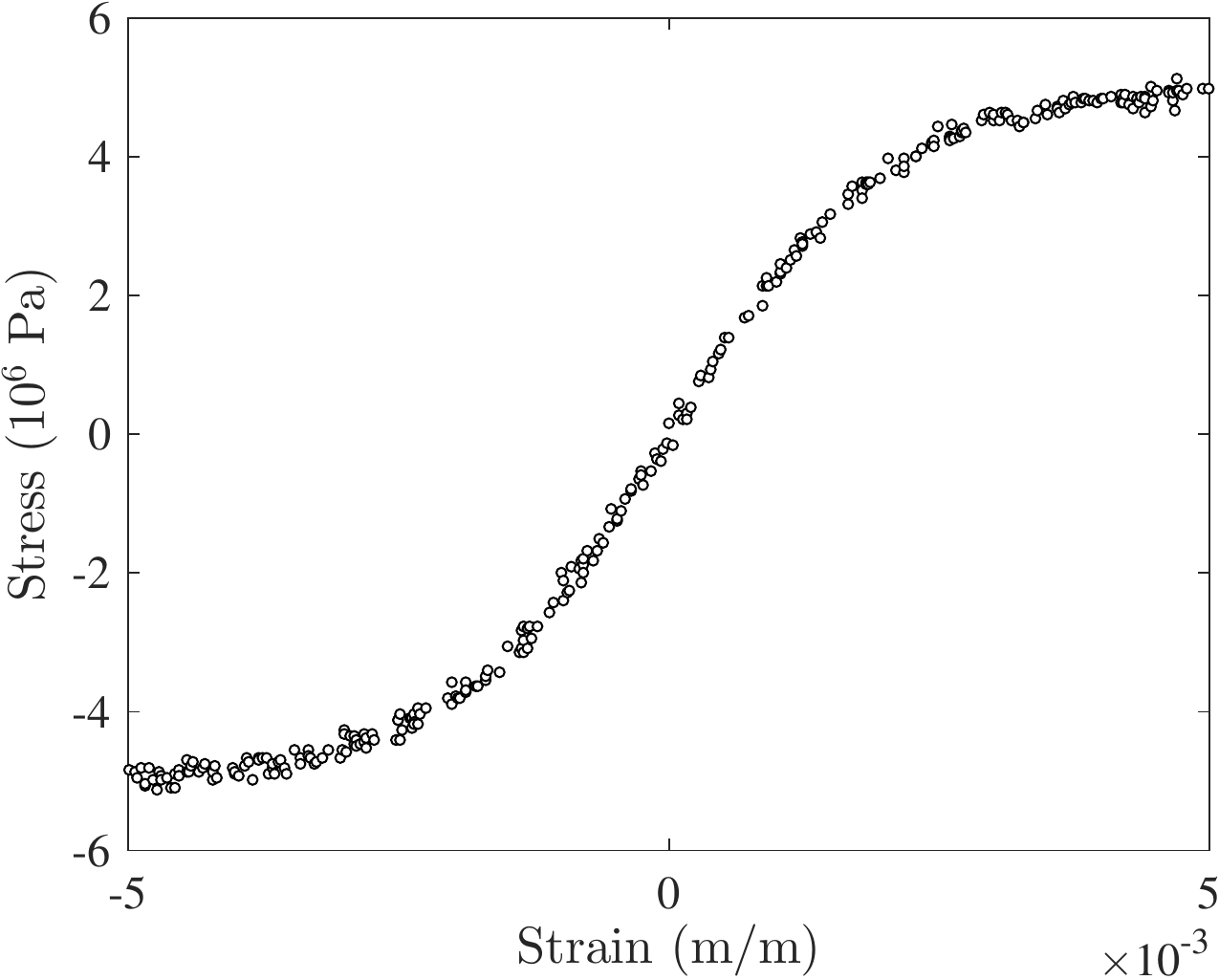}
  \caption{A material data set for the numerical experiments. }
  \label{fig:ten_bar_data_set}
\end{figure}

\begin{figure}[tp]
  \centering
  \includegraphics[scale=0.45]{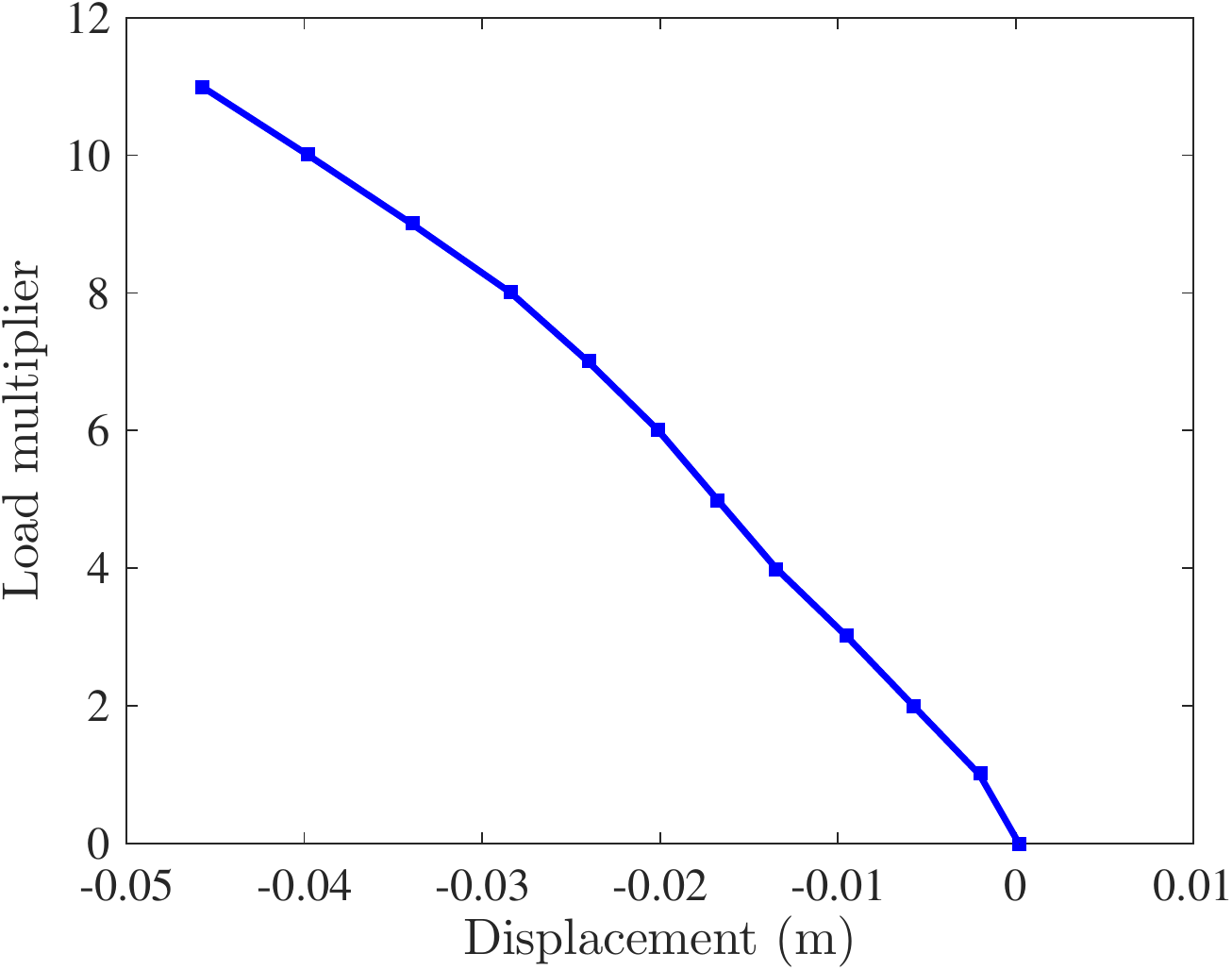}
  \caption{The obtained equilibrium path. }
  \label{fig:ten_bar_equilibrium_path}
\end{figure}

\begin{figure}[tp]
  \centering
  \begin{subfigure}[b]{0.47\textwidth}
    \centering
    \includegraphics[scale=0.50]{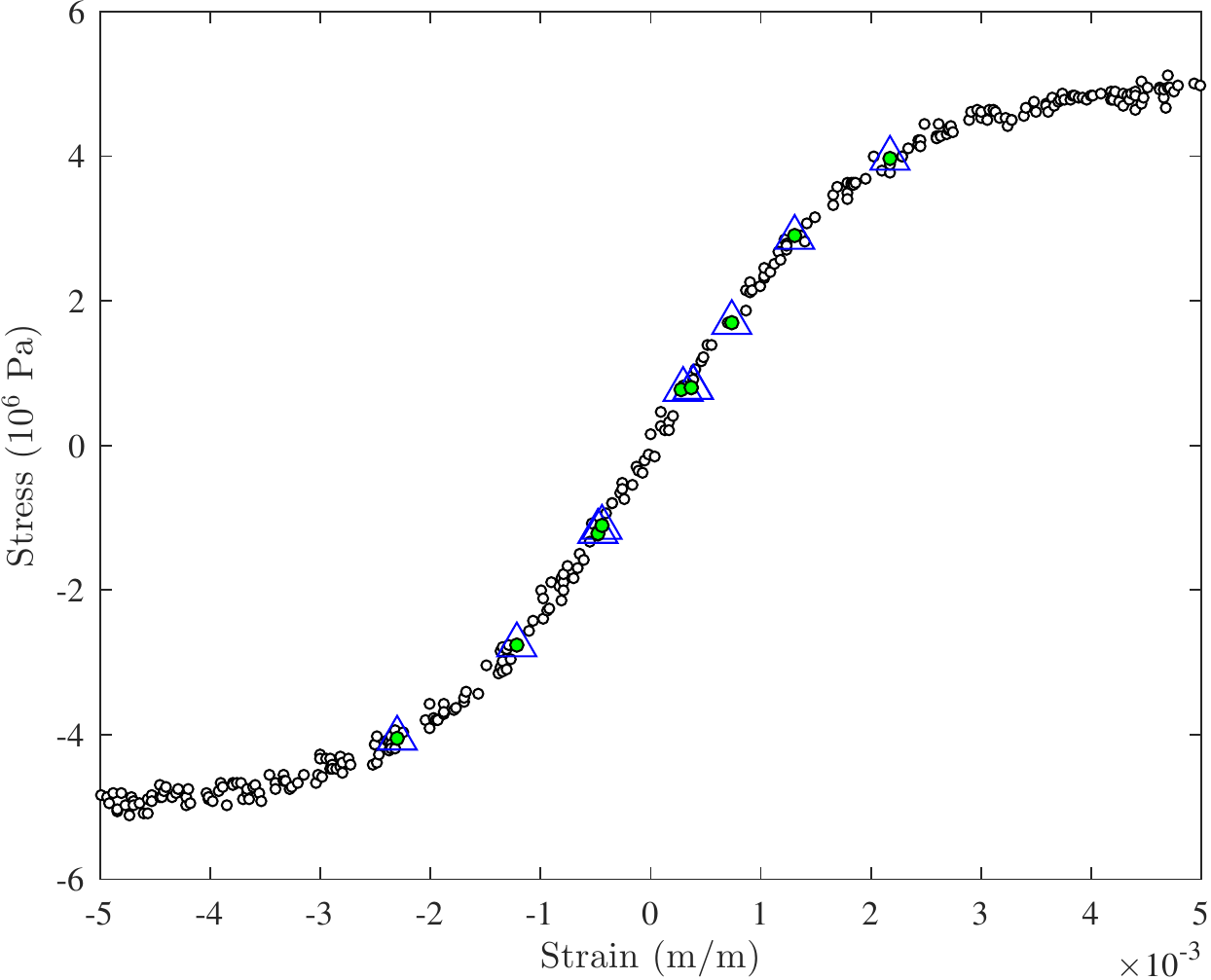}
    \caption{}
    \label{fig:ten_bar_load_11}
  \end{subfigure}
  \hfill
  \begin{subfigure}[b]{0.47\textwidth}
    \centering
    \includegraphics[scale=0.50]{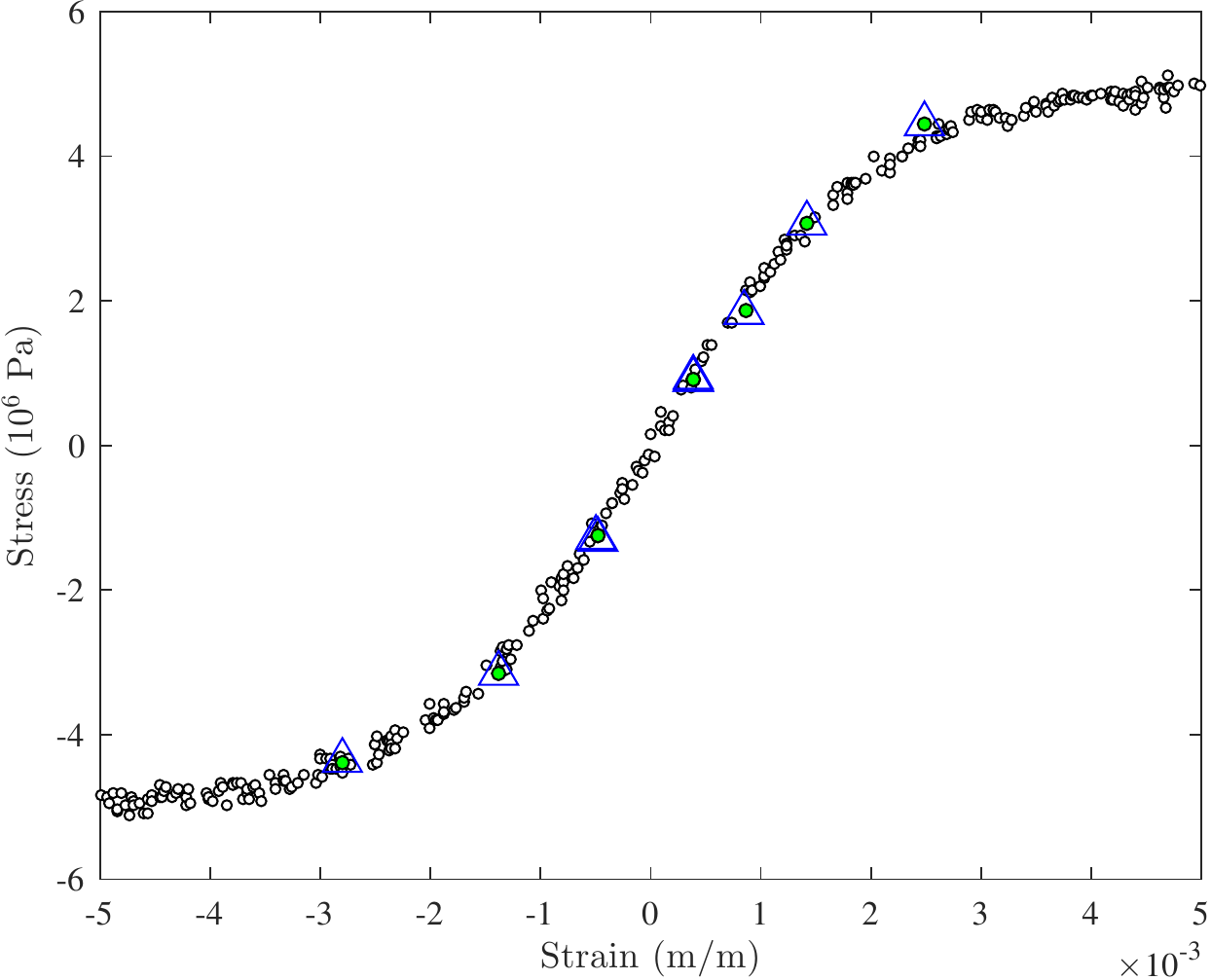}
    \caption{}
    \label{fig:ten_bar_load_12}
  \end{subfigure}
  \caption{The solutions obtained by the MIQP approach for 
  \subref{fig:ten_bar_load_11} $\lambda=10.0$; and 
  \subref{fig:ten_bar_load_12} $\lambda=11.0$. 
  ``{\em triangle\/}'' The stress and strain of each member; and 
  ``{\em filled circle\/}'' the nearest data points. 
  }
  \label{fig:ten_bar_load}
\par\bigskip
  \centering
  \begin{subfigure}[b]{0.47\textwidth}
    \centering
    \includegraphics[scale=0.50]{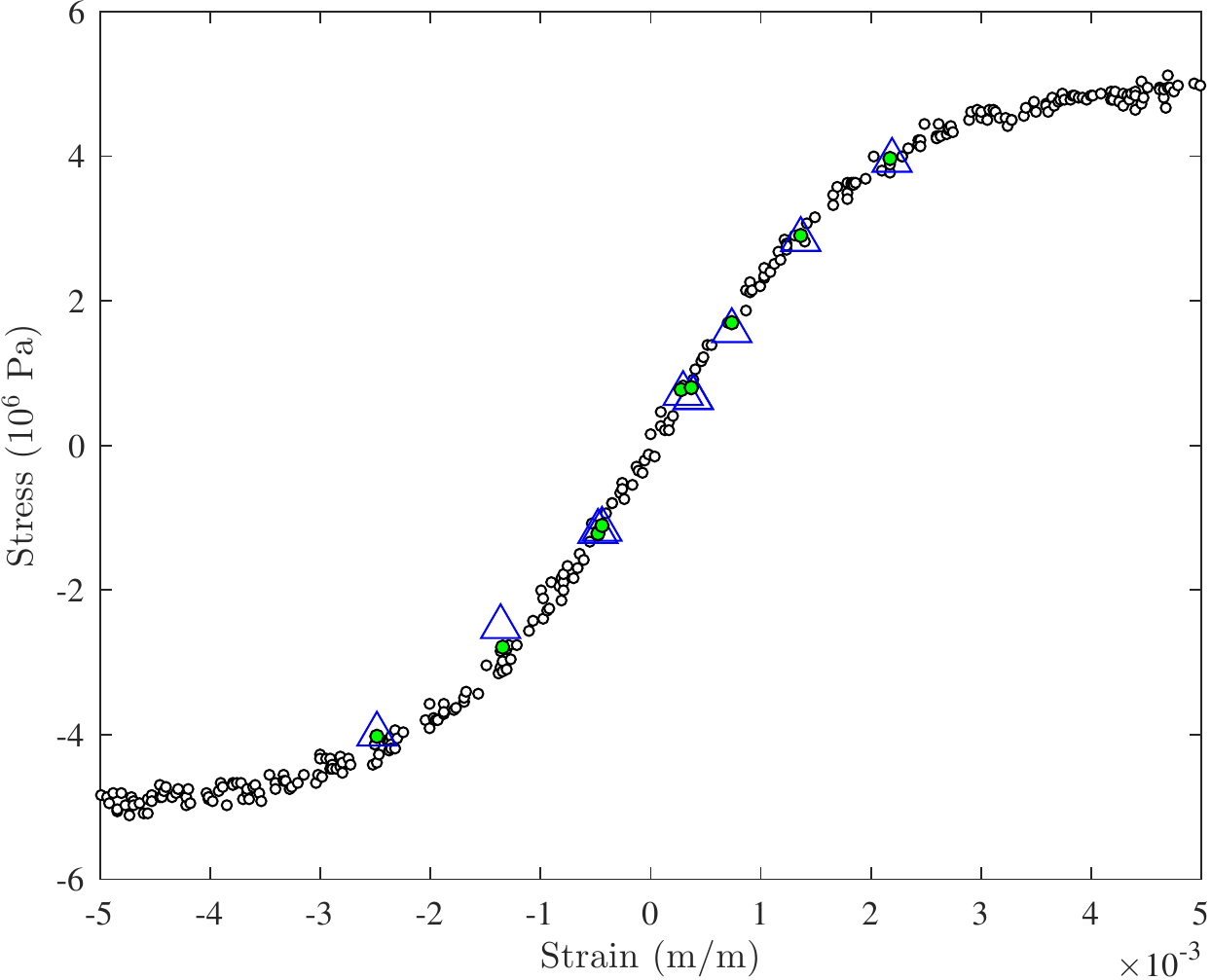}
    \caption{}
    \label{fig:ten_bar_heuristic_load_11}
  \end{subfigure}
  \hfill
  \begin{subfigure}[b]{0.47\textwidth}
    \centering
    \includegraphics[scale=0.50]{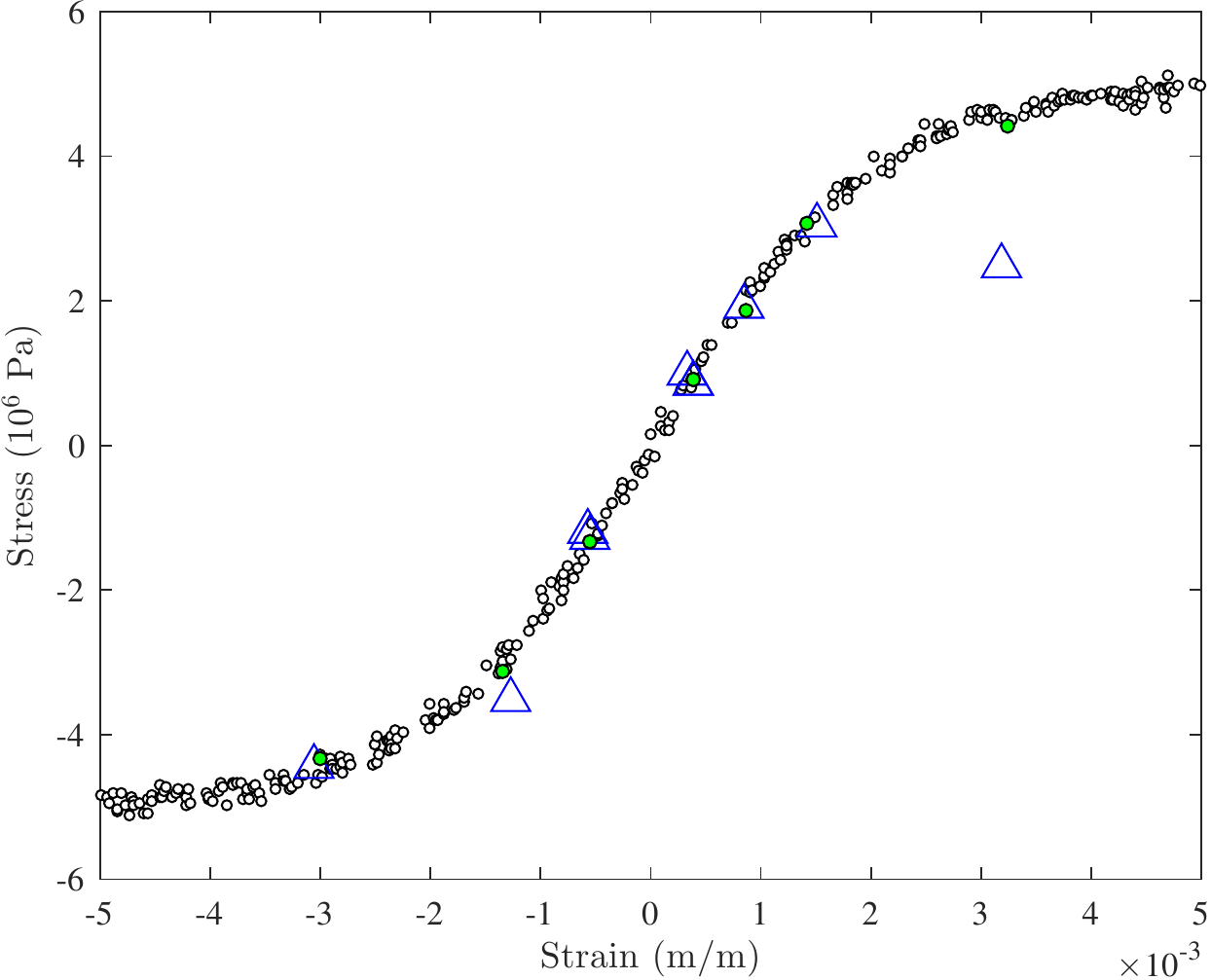}
    \caption{}
    \label{fig:ten_bar_heuristic_load_12}
  \end{subfigure}
  \caption{The solutions obtained by the heuristic in \cite{KO16} for 
  \subref{fig:ten_bar_load_11} $\lambda=10.0$; and 
  \subref{fig:ten_bar_load_12} $\lambda=11.0$. 
  ``{\em triangle\/}'' The stress and strain of each member; and 
  ``{\em filled circle\/}'' the nearest data points. 
  }
  \label{fig:ten_bar_heuristic_load}
\end{figure}

\begin{table}[bp]
  \centering
  \caption{Computational results.}
  \label{tab:result}
  \begin{tabular}{rrrrrr}
    \toprule
    & \multicolumn{3}{c}{MIQP} & \multicolumn{2}{c}{Heuristic \cite{KO16}} \\
    \cmidrule(lr){2-4}  \cmidrule(l){5-6}
    $\lambda$ & Opt.\ ($10^{-3}\,\mathrm{J}$) & Time (s) & {\#}BnB-node 
    & Obj.\ ($10^{-3}\,\mathrm{J}$) & {\#}iter.\ \\
    \midrule
    $0.0$  &  $6.528$ & $39.6$ & $10383$ & $1664.071$ & $10$ \\
    $1.0$  & $23.592$ & $29.4$ & $9321$  & $321.427$ & $17$ \\
    $2.0$  & $49.087$ & $29.6$ & $12417$ & $453.200$ & $235$ \\
    $3.0$  & $26.810$ & $48.6$ & $17647$ & $221.352$ & $66$ \\
    $4.0$  & $41.118$ & $87.2$ & $74230$ & $377.560$ & $46$ \\
    $5.0$  & $30.580$ & $33.4$ & $10589$ & $271.801$ & $338$ \\
    $6.0$  & $13.360$ & $33.6$ & $21677$ & $490.635$ & $27$ \\
    $7.0$  & $73.674$ & $62.5$ & $60549$ & $537.069$ & $619$ \\
    $8.0$  & $106.175$ & $43.5$ & $31200$ & --- & $(>10000)$ \\
    $9.0$  & $26.902$ & $36.8$ & $26276$ & $395.246$ & $25$ \\
    $10.0$ & $2.878$ & $63.3$ & $47830$  & $426.654$ & $22$ \\
    $11.0$ & $20.120$ & $52.5$ & $49763$ & $9011.320$ & $464$ \\
    \bottomrule
  \end{tabular}
\end{table}

As preliminary numerical experiments, the presented MIQP problem was 
solved with CPLEX ver.~12.8.0 \cite{cplex}. 
Computation was carried out on a $2.2\,\mathrm{GHz}$ Intel Core i5-5200 
processor with $8\,\mathrm{GB}$ RAM. 
We allowed CPLEX to use up to four threads. 
The integrality tolerance and 
the relative MIP gap tolerance (i.e., the tolerance on the relative gap 
between the objective value of the best feasible solution and that of 
the best branch-and-bound node remaining) of CPLEX were set to $0$, 
For comparison, Algorithm~1 in \cite{KO16} was also implemented 
in Matlab ver.\ 9.0. 
We set the initial point for this algorithm to the zero vector. 

Consider the planar truss shown in \reffig{fig:10bar}. 
This truss has $m=10$ members and $n=8$ degrees of freedom 
of the nodal displacements. 
As for the external load, 
vertical downward forces of $0.4 \lambda$ in $\mathrm{kN}$ are applied at the bottom two nodes as shown in \reffig{fig:10bar}, 
where $\lambda$ is the load multiplier. 
\reffig{fig:ten_bar_data_set} shows a material data set, which 
consists of $d=300$ data points. 
Hence, our MIQP problem in \eqref{P.original.1} has $m d=3000$ binary variables. 
We set constant $c$ in the objective function in \eqref{P.original.1.1} 
to the mean of 
$\check{\sigma}_{1}/\check{\varepsilon}_{1},\dots,\check{\sigma}_{d}/\check{\varepsilon}_{d}$, 
which yields $c=1.622\,\mathrm{GPa}$. 

The proposed method computed the equilibrium path shown in 
\reffig{fig:ten_bar_equilibrium_path}, which depicts 
the variation of the vertical displacement of the bottom rightmost node 
with respect to the load multiplier, $\lambda$.  
\reftab{tab:result} reports the computational results. 
For the MIQP approach, 
``opt.''\ in \reftab{tab:result} reports the optimal value found by CPLEX, 
``time'' is the computational time, and 
``{\#}BnB-node'' is the number of enumeration nodes explored by CPLEX. 
Also, for the heuristic in \cite{KO16}, 
 ``obj.''\ is the objective value of the solution found by the method, and 
``{\#}iter.''\ is the number of iterations. 
The computational time required by the MIQP approach is about or less 
than $60\,\mathrm{s}$. 
The heuristic \cite{KO16} does not necessarily converge; namely, 
for the problem instance with $\lambda=8.0$, it did not converge within 
$10000$ iterations. 
Also, even if it converges, the obtained solution is not guaranteed to 
be optimal; indeed, in all the converged cases, its solution has 
an objective value greater than the one computed by CPLEX. 

\reffig{fig:ten_bar_load} depicts typical solutions obtained by the 
proposed method. 
Here, a triangle indicates a pair of the member strain and stress, 
$(\varepsilon_{i},\sigma_{i})$, 
and a filled circle indicates its nearest data point, 
$(e_{i},s_{i})$. 
For the same problem instances, \reffig{fig:ten_bar_heuristic_load} 
shows the solutions obtained by the heuristic \cite{KO16}. 
It is worth noting that the objective values of these solutions are more 
than $100$ times larger than the optimal values.

\section{Conclusions}
\label{sec:conclude}

In this note, we have seen that a problem solved for 
data-driven computational mechanics \cite{KO16} can be formulated as 
a mixed-integer quadratic programming (MIQP) problem. 
We can solve the MIQP problem globally by using a standard 
optimization software package. 
Therefore, the results obtained by the presented approach can be 
utilized for benchmarking the other algorithms. 

\paragraph{Acknowledgments}

This work is partially supported by 
JSPS KAKENHI 17K06633 and 18K18898.

\end{document}